\documentclass{lmcs}
\pdfoutput=1

\usepackage{lastpage}
\lmcsdoi{22}{2}{2}
\lmcsheading{}{\pageref{LastPage}}{}{}%
{Sep.~25,~2023}{Apr.~07,~2026}{}

\usepackage[utf8]{inputenc}
\usepackage[T1]{fontenc}
\usepackage[british]{babel}
\usepackage{microtype}

\usepackage{amssymb,bbm,mathrsfs,stmaryrd}

\usepackage[pdfencoding=auto,pdfusetitle]{hyperref}

\usepackage{url}
\usepackage[shortlabels]{enumitem}

\usepackage{algorithm}
\usepackage[noend]{algpseudocode}

\usepackage{tikz}
\usetikzlibrary{matrix}
\usetikzlibrary{positioning}
\usetikzlibrary{calc}
\usetikzlibrary{arrows}
\tikzset{> =stealth}

\usepackage{etoolbox}
\newcommand{\addQEDstyle}[2]{\AtBeginEnvironment{#1}{\pushQED{\qed}\renewcommand{\qedsymbol}{#2}}\AtEndEnvironment{#1}{\popQED}}

\addQEDstyle{defi}{$\triangle$}
\addQEDstyle{exa}{$\triangle$}

\newenvironment{continueexample}[1]{\def\continuation{\ref*{#1}}
\newtheorem*{excont}{Example \continuation}
\addQEDstyle{excont}{$\triangle$}
\begin{excont}[Continued]}{\end{excont}}

\theoremstyle{definition}
\newtheorem{model}[thm]{Model}
\addQEDstyle{model}{$\triangle$}

\renewcommand{\epsilon}{\varepsilon}
\renewcommand{\phi}{\varphi}
\newcommand{\N}{\mathbb{N}}

\newcommand{\Q}{\mathbb{Q}} 
\newcommand{\R}{\mathbb{R}}

\renewcommand{\O}{\mathcal{O}}

\newcommand{\op}{^\mathrm{op}}

\newcommand{\inv}{^{-1}}
\newcommand{\id}{\mathrm{id}}
\newcommand{\powerset}{\mathcal{P}}
\newcommand{\powerfin}{\mathcal{P}_{\mathsf{fin}}}

\newcommand*{\dirsup}[1][\hspace{0.75ex}]{\operatorname*{\bigvee{}^{\!\uparrow}_{\hspace{-0.2ex} #1}}}

\newcommand{\Loc}{\mathrm{Loc}}
\newcommand{\Top}{\mathrm{Top}}
\newcommand{\Set}{\mathrm{Set}}
\newcommand{\Pos}{\mathrm{Pos}}
\newcommand{\Hom}{\mathrm{Hom}}

\newcommand{\ev}{\mathrm{ev}}
\newcommand{\evweak}{\widetilde{\mathrm{ev}}}

\begin{document}

\title[Machine Space I]{Machine Space I: Weak exponentials\texorpdfstring{\\}{} and quantification over compact spaces}

\thanks{The second author acknowledges financial support from the Centre for Mathematics of the University of Coimbra (UIDB/00324/2020, funded by the Portuguese Government through FCT/MCTES)}

\author[P.~F.~Faul]{Peter F.~Faul\lmcsorcid{0000-0001-5241-2603}}
\author[G.~Manuell]{Graham Manuell\lmcsorcid{0000-0002-3556-5293}}

\address{Stellenbosch University, Stellenbosch, South Africa}
\email{peter@faul.io, graham@manuell.me}

\begin{abstract}
Topology may be interpreted as the study of verifiability, where opens correspond to semi-decidable properties. In this paper we make a distinction between verifiable properties themselves and processes which carry out the verification procedure. The former are simply opens, while we call the latter \emph{machines}.
Given a frame presentation $\O X = \langle G \mid R\rangle$ we construct a space of machines $\Sigma^{\Sigma^G}$ whose points are given by formal combinations of basic machines corresponding to generators in $G$. This comes equipped with an `evaluation' map making it a weak exponential with base $\Sigma$ and exponent $X$.

When it exists, the true exponential $\Sigma^X$ occurs as a retract of machine space. We argue this helps explain why some spaces are exponentiable and others not.
We then use machine space to study compactness by giving a purely topological version of Escardó's algorithm for universal quantification over compact spaces in finite time.
Finally, we relate our study of machine space to domain theory and domain embeddings.
\end{abstract}

\maketitle

\section{Introduction}\label{sec:intro}

The interpretation of topology in terms of \emph{verifiable properties} provides insight into many aspects of the field. However, if this perspective is taken completely seriously it prompts some tricky questions concerning the failure of exponentials of topological spaces. In this paper, we will address these questions by introducing the notion of \emph{machine space} and demonstrate its utility by extending Escardó's algorithm for universal quantification over compact spaces (see \cite{escardo2007infinite}) from data types to general spaces.

A property $P$ is said to be \emph{verifiable} if whenever an element satisfies $P$, this can be established by finite means. Notably absent are any constraints on elements which do not satisfy $P$ --- that is, you need not be able to verify this fact. 

Let $X$ be some collection of points and let us consider what properties the set $\O X$ of verifiable propositions about $X$ should have. Given two propositions $U$ and $V$ we can verify their conjunction $U \wedge V$ by simply verifying $U$ and $V$ in turn. Furthermore, the constantly true proposition is trivially verifiable and so $\O X$ is closed under finite conjunctions.
Note that this argument does not extend to infinite conjunctions, since checking infinitely many verifiable propositions in turn is not achievable via finite means. However, verifiable properties are closed under disjunctions of arbitrary cardinality:
to verify that a disjunction holds at some point, we need only verify that one disjunct holds, which does not require any infinite processes. Note in particular that the constantly false proposition is vacuously verifiable, since we need only verify propositions when they hold.

In summary, the logic of verifiability admits finite conjunctions and arbitrary disjunctions. Identifying each proposition $U$ with the set of points which satisfy it recovers the usual notion of a \emph{topological space} with verifiable properties corresponding to open sets, conjunctions interpreted as set intersection and disjunctions as set union. Thus, point-set topology provides one formalisation of the theory of verifiability. Indeed, this perspective even extends to \emph{continuous functions}.

Let $f\colon X \to Y$ be some function between spaces, and let $U$ be some verifiable property on $Y$. Supposing $f$ is physically realisable, we can verify if $f(x) \in U$ or equivalently if $x \in f\inv(U)$. Thus, we recover that preimages of open sets must be open. Note that from this point of view, all physically realisable functions must necessarily be continuous.

Concrete topological spaces can be fruitfully interpreted from the perspective of\linebreak verifiability.
\begin{exa}
Consider the case of the real numbers $\R$. We might imagine having devices that can measure a quantity $x$, each having some fixed precision $\epsilon > 0$ and giving rational outputs. More precisely, for a real number $x$, such a device will (nondeterministically) give some $q \in \Q$ such that $\lvert x - q\rvert < \epsilon$.
Moreover, let us assume that for any precision there exists a device that is at least that precise.

We can use these to semi-decide if $x > a$ for some given $a \in \R$, since if $x > a$ then there is always some $\epsilon$ such that  $x-\epsilon > a$ and so using a device with such a precision will verify that $x$ is indeed greater than $a$. By a similar argument, we can semi-decide if $x < b$ and hence if $x \in (a,b)$.
On the other hand, it is not always possible to check if $x \ge a$, as in the scenario where $x = a$ we need infinite precision to be sure that $x$ is not actually very slightly smaller than $a$ and any given device has only finite precision.
This agrees with the familiar topology on $\R$.
\end{exa}

This interpretation of topology is due to Smyth in \cite{smyth1983power,symth1993topology} and has been expanded on subsequently in \cite{vickers1996topology} and various other papers discussed below. 

Given this perspective, it is possible to derive interpretations of various other topological concepts in terms of verifiability. For instance, it is well known that a compact set $K$ corresponds to a set that can be universally quantified over --- that is, given a verifiable property $P$ it is possible to verify the property that all members of $K$ satisfy $P$. (Note that $K$ may in fact be infinite!)

In this paper we make use of the \emph{pointfree} approach which emphasises the verifiable properties $P$ over the collections of points which satisfy them.
The analogue of topological spaces from this perspective are called \emph{locales}, though we will often simply call them `spaces'.

\subsubsection*{Philosophical problems}
 The verifiability interpretation raises the following problem of a philosophical nature: do the verifiable properties on a space themselves have verifiable properties so that they can be arranged into a topological space? Naively the answer would appear to be `yes'. For instance, let $x \in X$ be a point in a space $(X,\O X)$. It appears that we should be able to verify whether a verifiable property in $\O X$ contains $x$ or not --- as for any verifiable property $U$, by definition we can verify if an element $x$ satisfies it. While this is a compelling argument, the mathematics tells a different story. What we are asking for is equivalent to asking for a certain exponential object to exist in the category of topological spaces, but it is well known that these exponentials do not always exist. How can we reconcile this failure with our philosophical intuition?

 A seemingly unrelated problem, which was actually the original motivation behind the paper, concerns compactness: if compact spaces allow for universal quantification of open predicates, is there a \emph{uniform} procedure for universal quantification that works for all compact spaces? Or do different spaces require ad-hoc procedures?

 The work of Escardó indicates that it is likely the former. In \cite{escardo2007infinite,escardo2008exhaustible} he provides such a uniform algorithm for decidable properties in the context of programming language semantics. While this is compelling, it is not the complete solution to our problem. For one, Escardó is mainly interested in decidable properties, instead of the more general verifiable properties. Moreover, these papers concern Kleene--Kreisel spaces instead of general topological spaces or locales.
 
 Both of these problems are concerned with verifiable properties themselves and in what follows we provide a treatment of them that resolves both problems.
 
\subsubsection*{Mathematical solutions}
It can be instructive to imagine that for each open in a space, there exists some machine or program which carries out the verification procedure. Each machine takes points of the space as input and will either halt (in finite time) or run forever, depending on whether the point belongs to the associated open or not. Here we see an example of the analogy between verifiability in topology and semi-decidability in computability theory \cite{smyth1983power,escardo2004synthetic}. This link is also exhibited by the compactness algorithm, which we discuss later.

We might be tempted to identify machines with opens, which can in turn be seen to correspond to continuous maps from $X$ into the Sierpiński space $\Sigma$.
Since we are thinking of these machines as `real things' we can interact with, we can expect there to be natural verifiable properties concerning the machines themselves (for example, does the machine halt on some given point?) yielding a \emph{space of machines} for each space $X$.
However, if $X$ is not locally compact, the exponential object $\Sigma^X$, which we would think of as the space of opens, does not exist and so we are forced to distinguish between machines and opens --- unlike functions from $X$ to $\Sigma$, machines are not extensional. Moreover, since spaces of machines always exist, machines are better than opens for some purposes.

In order to formalise these ideas, we describe an explicit construction of a space of machines. We replace the (potentially nonexistent) space of opens $\Sigma^X$ with a certain \emph{weak exponential}.
Our starting point is to fix a presentation for $X$ with generators $G$. (From the classical perspective these can be viewed as a set of subbasic opens $G$ for the space $X$.) Since $G$ is a set, the space $\Sigma^G$ exists. The space $X$ embeds canonically into $\Sigma^G$ and so the opens of $X$ are restrictions of opens in $\Sigma^G$. We will view distinct opens of $\Sigma^G$ which restrict to the same open $U$ in $X$ as distinct machines which accept (i.e.\ halt on) precisely the elements of $U$.
A crucial point is that the space $\Sigma^G$ is locally compact and we can take $\Sigma^{\Sigma^G}$ to be the space of machines of $X$. This space may be thought of as a reasonably canonical weak exponential associated to $X$ with base $\Sigma$. 
More concretely, the points of $\Sigma^{\Sigma^G}$ can be thought of as formal joins of formal meets of the generators, which are represented by certain programs that run the `basic machines' from $G$ in parallel.

We then relate machine space to $\Sigma^X$, when the latter exists. In particular, there is a canonical quotient map $\Sigma^{\Sigma^G} \twoheadrightarrow \Sigma^X$ sending each machine to its associated open. Under the interpretation so far described an open merely represents a verifiable property in an abstract sense, whereas a machine is some process that concretely semi-decides memberships of elements. Given an open and a point in a general space $X$ there is \emph{not} an obvious way to verify that the point lies in the open, as the open alone does not give any such procedure. This helps explain the fact that for general spaces the collection of opens equipped with the Scott topology does not have a continuous evaluation map. We could however expect the evaluation map to be continuous if there were some way to associate a machine to each open. Indeed, we show that when $X$ is locally compact the canonical quotient map of machine space onto $\Sigma^X$ always has a section, allowing a continuous assignment of opens to machines which represent them.

One way in which machine space is useful is in understanding compactness.
From the perspective of verifiability a \emph{compact space} is intuitively a space that can be universally quantified over in finite time \cite{taylor2011foundations,escardo2005notes}. More explicitly, if $P$ is some verifiable property and $K$ is some compact space then the question of whether all the members of $K$ satisfy $P$ is verifiable. That it is sometimes possible to universally quantify over infinite spaces is especially interesting and captures the intuitive idea that `compact spaces behave like finite sets'.
This idea can be formalised using hyperdoctrines \cite{lawvere1969adjointness,pittsLogic,manuellThesis}, but one might also consider an alternative approach involving the space of opens to have some appeal.

If $K$ is a compact, locally compact space, then we can check if a verifiable property holds on all of $K$ and so we should be able to verify whether the corresponding open in $\Sigma^K$ is equal to the largest element of $\Sigma^K$, namely $K$ itself. In other words, $K$ is compact if and only if the singleton $\{K\}$ is open in $\Sigma^K$ (see \cite{escardo2005notes}).
Of course, this viewpoint is inapplicable outside of the locally compact case, since the exponential $\Sigma^K$ will not exist \cite{hyland1981exponentials}.

As might be expected, this can be extended to the general setting of compact (but not necessarily locally compact) locales by replacing $\Sigma^K$ with machine space. This perspective on the universal quantification is essentially due to Escardó \cite{escardo2004synthetic}. Instead of asking if an open of $K$ equals $K$, we ask if an open in $\Sigma^G$ \emph{contains} $K$ --- that is, if the machine accepts all points of $K$. By the Hofmann--Mislove theorem \cite{hofmann1981local}, there is an associated open corresponding to all of the machines which cover $K$, which plays the same role as the singleton $\{K\}$ in the previous approach.

Our contribution is obtaining this open via an explicit algorithm for universally quantifying over a compact space. This is a very general and purely topological/localic version of Escardó's algorithm \cite{escardo2004synthetic,escardo2007infinite} for universal quantification over Cantor space.
(See also \cite{escardo2008exhaustible} which discusses universal quantification over data types in the setting of higher type computation.)
Unlike Escardó, we do not work within a particular formal semantics, but remain agnostic towards the particular `model' of the topological structures.

Finally, we discuss some links with domain theory, give some concrete examples and show how the original algorithm for quantification over Cantor space can be recovered from ours.

In a later paper we intend to explore the links between machine space and powerlocales.

\section{Background}

\subsection{Pointfree topology}

The logic of verifiability can also be studied purely through the verifiable propositions themselves, without any regard to the points which may or may not satisfy these propositions. This is the perspective of \emph{pointfree topology}, the study of topology through the lattice of open sets. In this section we introduce some of the basic notation and concepts which will be used throughout the rest of the paper. See \cite{picado2011frames,vickers1996topology} for further details.

We have seen that verifiable properties are closed under finite conjunctions and arbitrary disjunctions. Identifying logically equivalent propositions and ordering the equivalence classes by logical entailment we obtain a lattice. This is the Lindenbaum--Tarski algebra for the logic. For the conjunctions and disjunctions to behave as expected, we require that this lattice is distributive. This motivates the following definition.

\begin{defi}
 A \emph{frame} is a complete lattice which satisfies the distributivity law \[u \wedge \bigvee_{i \in I} v_i = \bigvee_{i \in I} (u \wedge v_i).\]
 A \emph{locale} $X$ is formally the same thing as a frame, but is thought of as being the abstract space that has the corresponding frame $\O X$ as its lattice of opens.
\end{defi}

Of course, the open sets of any topological space form a frame under intersection and union. Examples of this form are called \emph{spatial}. There are also non-spatial examples (such as complete atomless Boolean algebras).

We can also describe continuous functions from this perspective. Since the preimage of an open is an open, every continuous function induces a function between the lattices of open sets but in the opposite direction. From set-theoretic properties of preimages, we see that this map preserves finite meets and arbitrary joins. Thus, we define a locale morphism as follows.

\begin{defi}
 A \emph{frame homomorphism} $f\colon L \to M$ is a function $f\colon L \to M$ between frames that preserves finite meets and arbitrary joins. A \emph{locale morphism} $f\colon X \to Y$ is a frame homomorphism $f^*\colon \O Y \to \O X$ between the corresponding frames of opens.
\end{defi}

In good situations, a topological space can be recovered from its frame of open sets. (See \cite{picado2011frames} for more details.)
In particular, we can talk about \emph{points} of a locale.
As in $\Top$, we can identify points of a locale with maps from the terminal object $1$.
These correspond to frame homomorphisms from $\O X$ to $\O 1 \cong \{0,1\}$, which can be understood logically as assigning truth values to each verifiable proposition $P$ corresponding to whether $P$ holds at~$x$.

Frames are well-behaved (infinitary) algebraic structures and can be \emph{presented} by generators and relations (see \cite{Johnstone1982stone}).
In particular, free frames exist. The free frame on $G$ can be described explicitly as the frame of downsets on the meet-semilattice of finite subsets of $G$ ordered by reverse inclusion. The finite subsets of $G$ are interpreted as formal meets of generators, while the downsets are viewed as formal joins of these. Note that every formal expression involving frame operators on the generators can be brought into the form of a join of finite meets by distributivity.

\begin{defi}
 A \emph{presentation} for a frame consists of a set of generators $G$ and a set of relations $R$ consisting of formal equalities (or inequalities) between formal combinations of generators --- for example, $g_1 \wedge g_2 = \bigvee_{i = 3}^\infty g_i \wedge g_2$. Explicitly, $R$ can be defined to be an (arbitrary) subset of $F(G) \times F(G)$ where $F(G)$ is the free frame on the set $G$.
 
 The frame $\langle G \mid R\rangle$ defined by such a presentation contains elements corresponding to the generators $g \in G$ and which satisfy the relations given in $R$. Moreover, it is the initial frame satisfying this property: for any frame $M$ and function $f\colon G \to M$ for which the images of the generators under $f$ satisfy the relations from $R$, there is a unique frame homomorphism $\bar{f}\colon \langle G \mid R \rangle \to M$ making the triangle commute.
 \begin{center}
  \begin{tikzpicture}[node distance=3.5cm, auto]
    \node (A) {$G$};
    \node (B) [above of=A] {$\langle G \mid R\rangle$};
    \node (C) [right of=A] {$M$};
    
    \draw[->] (A) to node {} (B);
    \draw[->] (A) to node [swap] {$f$} (C);
    \draw[dashed,->] (B) to node {$\bar{f}$} (C);
  \end{tikzpicture}
\end{center}
Explicitly, the frame $\langle G \mid R\rangle$ can be constructed as the quotient of the free frame $F(G)$ by the congruence relation generated by $R$.

Note that taking $M$ in the diagram above to be the initial frame $\{0,1\}$, we see that the points of $\langle G \mid R\rangle$ are given by specifying the subset of the generators that are viewed as `true' such that relations become logical formulae. For example, our relation above becomes $g_1 \wedge g_2 \iff \exists i \ge 3.\ g_i \wedge g_2$, so $g_1$ and $g_2$ should be true if and only if $g_i$ and $g_2$ are true for some $g_i \ge 3$.
\end{defi}

A general presentation $\langle G \mid R\rangle$ can be obtained from a coequaliser $F(R) \rightrightarrows F(G) \twoheadrightarrow \langle G \mid R\rangle$ where the maps $F(R) \rightrightarrows F(G)$ are obtained by the universal property of $F(R)$ from the natural maps $R \rightrightarrows F(G)$.

One particularly important locale is $\Sigma$ corresponding to \emph{Sierpiński space}.
As a topological space this is $\{\bot,\top\}$ with the topology generated by the single open $\{\top\}$.
It can be formally defined as the locale corresponding to the free frame on a single generator (corresponding to $\{\top\}$).
Locale maps $u\colon X \to \Sigma$ correspond to opens of $X$ by the universal property of the free frame. Spatially, this can be understood as taking the preimage of $\{\top\}$. From the verifiability perspective this allows us to interpret
the elements of $\Sigma$ as corresponding to whether the verification process halts ($\top$), or runs forever ($\bot$).

Note that if $G$ is a set then $\Sigma^G$ represents the product $\Pi_{g \in G} \Sigma$. (We can alternatively describe it as an exponential object as defined in the following subsection, with $G$ viewed as a discrete locale.) The corresponding frame $\O(\Sigma^G)$ is precisely $F(G)$, the free frame on $G$. So if $X$ is a locale with $\O X = \langle G \mid R\rangle$, we can dualise the coequaliser for $\langle G \mid R \rangle$ to obtain the locale equaliser
\[X \hookrightarrow \Sigma^G \rightrightarrows \Sigma^R.\]

The equaliser morphism from $X \hookrightarrow \Sigma^G$ can be viewed as analogous to an embedding of topological spaces. Regular subobjects of locales are called \emph{sublocales} and correspond to frame quotients.

Another concept we recall here is that of compactness. The usual topological definition in terms of finite subcovers works equally well in this setting.

\begin{defi}
A locale $X$ is compact if whenever $\bigvee_{i \in I} u_i = 1$ there exists a finite set $F \subseteq I$ such that $\bigvee_{i \in F} u_i = 1$.
\end{defi}

As mentioned before, compactness has a very nice interpretation in terms of universal quantification of verifiable properties. We will explore this topic further in \autoref{sec:compactness}.

\subsection{The Scott topology and exponentials}

The set of opens $\O X$ of a locale $X$ is itself often endowed with the \emph{Scott topology}.
\begin{defi}
 A subset $V \subseteq \O X$ is \emph{Scott-open} if it is upward closed and if whenever a join $\bigvee D$ of a directed set $D$ is in $V$, then some $d \in D$ lies in $V$. Recall that a set $D$ is \emph{directed} if every finite subset of $D$ has an upper bound in $D$.
\end{defi}

The Scott topology can actually be defined more generally on \emph{dcpos}, which are often used to model data types in programming languages (as discussed in \autoref{sec:dcpo}).
\begin{defi}
 A poset which admits joins of all directed subsets is called a \emph{dcpo}. We will write $\dirsup S$ for the join of a directed set $S$.
 A morphism of dcpos or a \emph{Scott-continuous function} is a monotone map which preserves directed suprema.
\end{defi}

The order on a dcpo can be thought of as a definability order: the least element (if it exists) corresponds to a purely divergent computation, while maximal elements are completely defined. For more on this topic see \cite{abramsky1991domain,abramsky1994domain, cartwright2016domain,goubault2013nonhausdorff}.

\begin{defi}\label{def:continuous_dcpo}
The \emph{way-below relation} $\ll$ on a dcpo $L$ is defined so that $a \ll b$ if and only if whenever $b \le \dirsup D$ for a directed set $D$, there is some $d \in D$ with $a \le d$.
We say $L$ is \emph{continuous} if for every $b \in L$ we have $b = \dirsup \{a \in L \mid a \ll b\}$. (Note, in particular, that the set $\{a \in L \mid a \ll b\}$ is required to be directed).

A locale $X$ is \emph{locally compact} if and only if its frame of opens $\O X$ is continuous.
\end{defi}

Recall if $X,B$ are objects in a category that the \emph{exponential object} $B^X$ is equipped with a morphism $\ev\colon B^X \times X \to B$ such that for every $f\colon Y \times X \to B$ there is a unique $\overline{f}\colon Y \to B^X$ such that the following diagram commutes.
\begin{center}
  \begin{tikzpicture}[node distance=3.5cm, auto]
    \node (A) {$Y \times X$};
    \node (B) [below of=A] {$B^X \times X$};
    \node (C) [right of=B] {$B$};
    \draw[->] (A) to node [swap] {$\overline{f} \times X$} (B);
    \draw[->] (B) to node [swap] {$\ev$} (C);
    \draw[->] (A) to node {$f$} (C);
  \end{tikzpicture}
\end{center}
The exponential is thought of as the object of functions from $X$ to $B$ and $\ev$ as the evaluation map. In the diagram above $\overline{f}$ can be viewed as the curried version of $f$. We say $B$ is \emph{exponentiable} if the exponential $B^X$ exists for all objects $X$. 

If $X$ is a locally compact locale, then $X$ is exponentiable (see \cite{hyland1981exponentials}). In particular, cale $\Sigma^X$ exists and its frame is given by the frame of opens of the Scott topology on $\O X \cong \Hom(X, \Sigma)$. Furthermore, $\Sigma^X$ is itself locally compact (for instance, see \cite[Theorem 11.1]{vickers2004double}).
However, if $X$ is not locally compact, the exponential $\Sigma^X$ does \emph{not} exist. If $X$ is a non-locally-compact topological space, we find that the evaluation map $\ev\colon \O X \times X \to \Sigma$, which sends $(f,x)$ to $f(x)$, is not continuous when $\O X$ is given the Scott topology.

\section{Machine space}

We first tackle the problem of exponentials. We construct a `canonical weak exponential' associated to a frame presentation and relate these to true exponentials when they exist. Along the way, we discuss the philosophical reasons behind why weak exponentials with base $\Sigma$ exist, but true exponentials might not.

\subsection{Presentations and weak exponentials}

For a general space $X$ the `space of opens' $\Sigma^X$ need not exist. However, given a presentation $\O X = \langle G \mid R \rangle$ we will see that there is always a natural weak exponential given by $\Sigma^{\Sigma^G}$.

Recall that $\langle G \mid R \rangle$ is a quotient of the free frame on $G$, which in turn is the frame of opens of $\Sigma^G$.
The frame quotient $\Sigma^G \twoheadrightarrow \langle G \mid R \rangle$ exhibits $X$ as a sublocale of $\Sigma^G$. The points of $\Sigma^G$ correspond to subsets of $G$ and we can then interpret this inclusion as sending points of $X$ to the set of generating opens in which they lie.
The points in the image of this inclusion correspond to the frame homomorphisms $\O(\Sigma^G) \to \O 1$ which respect the relations. The additional points in $\Sigma^G$ may be thought of as `imaginary points' in the same sense that the exponential spaces used in (the different approach of) \cite[Chapter 10]{escardo2004synthetic} are described as imaginary exponentials. Just as complex numbers may be used to more easily prove results about real numbers, so will these imaginary points be useful for our study of compactness in \autoref{sec:compactness}.

The space $\Sigma^G$ is always locally compact and so we can form $\Sigma^{\Sigma^G}$ --- the space of opens of $\Sigma^G$. 
Of course, $\Sigma^{\Sigma^G}$ is a spatial locale whose points are Scott-opens of $\O X$ and whose topology is again given by the Scott topology. Alternatively, we can view $\Sigma^{\Sigma^G}$ as the double powerlocale on $G$ (see \cite{vickers2004double,vickers2004universal}), the topology of which is generated by opens ${\boxtimes} U$ for each $U \subseteq G$ where an element $m = \bigvee_{i \in I} \bigwedge_{j\in J_i} g_j$ of $\Sigma^{\Sigma^G}$ lies in ${\boxtimes} U$ if and only if there exists some $J_k$ such that $g_j \in U$ for each $j \in J_k$.

We now show that $\Sigma^{\Sigma^G}$ can be made into a \emph{weak exponential}. (Compare the construction of weak exponentials in $\Top$ given in \cite{rosicky1999}.)

\begin{defi}\label{def:weak_exponential}
 Let $X$ and $B$ be objects in some category. We say $Z$ is a \emph{weak exponential} with base $B$ and exponent $X$ if there exists a morphism $\evweak\colon Z \times X \to B$ such that for any $f\colon Y \times X \to B$ there exists a \emph{not necessarily unique} $\overline{f}\colon Y \to Z$ such that the following diagram commutes.
 \begin{center}
  \begin{tikzpicture}[node distance=3.5cm, auto]
    \node (A) {$Y \times X$};
    \node (B) [below of=A] {$Z \times X$};
    \node (C) [right of=B] {$B$};
    \draw[->] (A) to node [swap] {$\overline{f} \times X$} (B);
    \draw[->] (B) to node [swap] {$\evweak$} (C);
    \draw[->] (A) to node {$f$} (C);
  \end{tikzpicture}
\end{center}
In this case, we call $\evweak$ the \emph{evaluation map} for the weak exponential.
\end{defi}

\begin{prop}\label{prop:weak_exponential}
Let $X$ have presentation $\langle G \mid R\rangle$. The space $\Sigma^{\Sigma^G}$ together with the evaluation map $\evweak\colon \Sigma^{\Sigma^G} \times X \to \Sigma$ given by the composite $\Sigma^{\Sigma^G} \times X \hookrightarrow \Sigma^{\Sigma^G} \times \Sigma^G \xrightarrow{\ev} \Sigma$ is a weak exponential with base $\Sigma$ and exponent $X$.
\end{prop}
\begin{proof}
Since $X$ embeds into $\Sigma^G$ via $i_X\colon X \hookrightarrow \Sigma^G$, the product $\Sigma^{\Sigma^G} \times X$ embeds into $\Sigma^{\Sigma^G} \times \Sigma^G$. But $\Sigma^G$ is locally compact and so we have an evaluation map $\ev\colon  \Sigma^{\Sigma^G} \times \Sigma^G \to \Sigma$.
To see that $\evweak = \ev \circ (\Sigma^{\Sigma^G} \times i_X)$ makes $\Sigma^{\Sigma^G}$ a weak exponential we consider a morphism $u\colon A \times X \to \Sigma$ and must construct a morphism $v$ such that $u = \evweak \circ (v \times X)$. Consider the following diagram.

\begin{center}
  \begin{tikzpicture}[node distance=4cm, auto]
    \node (A) {$A \times X$};
    \node (A') [right of=A] {};
    \node (B) [below=2cm of A'] {$A \times \Sigma^G$};
    \node (B') [right of=B] {};
    \node (C) [below=2cm of B'] {$\Sigma$};
    \node (D) [left of=C] {$\Sigma^{\Sigma^G} \times \Sigma^G$};
    \node (E) [left of=D] {$\Sigma^{\Sigma^G} \times X$};
    \draw[right hook->] (A) to node {$A \times i_X$} (B);
    \draw[dashed,->] (B) to node {$u'$} (C);
    \draw[->] (A) to node [swap] {$\overline{u}' \times X$} (E);
    \draw[->] (B) to node [swap] {$\overline{u}' \times \Sigma^G$} (D);
    \draw[right hook->] (E) to node [swap] {$\Sigma^{\Sigma^G} \times i_X$} (D);
    \draw[->] (D) to node [swap] {$\ev$} (C);
    \draw[out=0,in=100,->] (A) to node {$u$} (C);
  \end{tikzpicture}
\end{center}

Observe that $u$ defines an open in $A \times X$. Since $A \times X$ is a subspace of $A \times \Sigma^G$, there exists a (not necessarily unique) open $u'$ in $A \times \Sigma^G$ which restricts to give $u$. By the universal property of the exponential $\Sigma^{\Sigma^G}$, we obtain a morphism $\overline{u}'$ making the bottom right triangle commute. The left-hand trapezium commutes since both composites give $\overline{u}' \times i_X$ and then a simple diagram chase confirms that $u = \evweak \circ (\overline{u}' \times X)$.
\end{proof}

The above construction is summarised in the following definition.
\begin{defi}
Let $X$ be a locale with $\O X = \langle G \mid R \rangle$. We define the associated \emph{(idealised) machine space} to be $\Sigma^{\Sigma^G}$ together with the evaluation map $\evweak\colon \Sigma^{\Sigma^G} \times X \to \Sigma$ defined in \autoref{prop:weak_exponential}.
\end{defi}
We explore the intuition behind the concept of machine space in the following section.

\subsection{Philosophical Interlude}\label{sec:philo}

In this section we examine the relationship between opens, i.e.\ verifiable properties, and the \emph{processes} that semi-decide them. We view these processes as being carried out by \emph{machines}, which we will identify with elements of the weak exponential $\Sigma^{\Sigma^G}$.

Given a space $X$ with $\O X = \langle G \mid R \rangle$ we interpret each generator $g \in G$ as a basic black-box machine which takes as input the points of $X$ and either halts after some finite amount of time, or runs forever. More complicated machines can be constructed from these basic machines. For instance, we may construct a composite machine by taking a finite collection of machines and insist that each machine halts on the given point before the composite machine is considered to have halted. We can also construct a further machine by running any number of these composite machines in parallel and halting on a point if at least one branch halts. 
This corresponds to taking formal joins of finite formal meets of generators and
so we see that these machines are precisely the elements of the free frame on $G$ -- that is, the points of the machine space $\Sigma^{\Sigma^G}$.

\begin{defi}\label{def:machine}
Given a space $X$ with $\O X = \langle G \mid R \rangle$, an \emph{(idealised) machine} $m$ over $X$ is a formal join of finite formal meets of generators written $m = \bigvee_{i \in I} \bigwedge_{j\in J_i} g_j$ with each $J_i$ finite.
\end{defi}

Notice that $\Sigma^{\Sigma^G}$ by itself contains no information about the behaviour of the machines on inputs from the locale $X$ with presentation $\langle G \mid R\rangle$, as this data is contained in the relation $R$. We can recover this information from the evaluation map $\evweak\colon \Sigma^{\Sigma^G} \times X \to \Sigma$ of \autoref{prop:weak_exponential}, which records whether a machine halts on a given point.

\begin{rem}
 In fact, for many purposes \emph{any} weak exponential with base $\Sigma$ and exponent $X$ behaves like a machine space for $X$. This could be interesting to explore further. However, for this paper we restrict our consideration to $\Sigma^{\Sigma^G}$, since the compactness algorithm given in \autoref{sec:compactness} is dependent on the presentation of the locale by generators and relations.
\end{rem}

A point $x \in X$ is accepted by a machine $m$ if there exists some $J_k$ such that $x$ is accepted by each `basic machine' $g_j$ for all $j \in J_k$ (i.e.\ such that $x$ lies in the open given by $g_j$ for all such $j$). In this way, each machine $m$ will define an open on which it halts.
Many of these machines will correspond to the same opens, as described by the relation $R$, so that each open in $X$ corresponds to an equivalence class of machines. In \autoref{def:quotient} of \autoref{sec:machine_spaces_and_space_of_opens} we will see that when the exponential $\Sigma^X$ exists there is a quotient map $\Sigma^{\Sigma^G} \twoheadrightarrow \Sigma^X$ assigning machines to opens.

\begin{rem}\label{rem:axiomatic_vs_verifiable}
Note that such relations must be taken as given axiomatically and cannot be verified to hold. (After all, we cannot show in finite time that two machines halt on precisely the same points.)
It is important to maintain a sharp distinction between the properties we assume as axioms, which define the spaces under consideration and cannot be proven empirically, and the properties we observe about the spaces we describe, which concern unspecified points and must be verifiable.
\end{rem}

Recall that the failure of the exponential $\Sigma^X$ to always exist poses a puzzle for the verifiable interpretation of topology. Luckily, we have seen that machine spaces always exist. It is now helpful to consider in detail how such spaces might correspond to concrete objects that can be learnt about via physical measurements. Note that with such ``real-world'' models the set of generators is assumed to be countable.

\begin{model}\label{mod:bot}
Let $X$ be a space with $\O X = \langle G \mid R \rangle$. We take each generator $g \in G$ to be a black box which accepts points of $X$ as input and either halts after some finite time or runs forever. We imagine these generators arranged in a line. Additionally there are a number of mobile machines which themselves accept points of $X$ as input and will move from generator to generator testing the point in the following manner.

The machine progresses through different stages. On the $i^\text{th}$ stage it selects a group of finitely many generators, writes the number $i$ onto them and (if they are not already running) runs them with the point as input. The generators are left to run as the next stage is started. If at any point every generator in a group has halted, the machine halts.

Such a machine instantiates the element $\bigvee_{i \in I} \bigwedge_{j\in J_i} g_j$ where each $J_i$ is the finite set of generators selected at stage $i$.

To accommodate the fact that $\vee$ and $\wedge$ are commutative, absorptive and idempotent we require the following of our machines.
\begin{enumerate}
    \item The order in which the machine visits each stage and the generators therein is chosen non-deterministically,
    \item If it visits the generators $S \subseteq G$ in a stage, then in some other stage it must visit $T \supseteq S$, for all such finite $T \subseteq G$,
    \item Finally, no two stages will have the exact same generators, and no stages will contain the same generator twice.
\end{enumerate}
 Note that associativity and distributivity are automatically satisfied with distributivity following from the fact that the formal operations are already reduced to the canonical form of a join of meets. These axioms together are then enough to give a frame.
\end{model}

What properties can we verify about the machines in this model? Of course, for a real implementation there are a number of irrelevant properties, such as the precise position of the machines, but taking an idealised perspective and viewing non-deterministic choices as conferring no additional information, there are still some nontrivial verifiable properties. This indicates that these machines are not complete black boxes and that there is some information that leaks while they run:
we can semi-decide if there exists some parallel branch of a machine $m$ such that each of the generators that run in that branch are contained in $U$. This is exactly the open ${\boxtimes} U$. Recall that the opens of this form generate the frame $\O \Sigma^{\Sigma^G}$.

We can now answer the question posed in the introduction. While we have seen that machine space is reasonably concrete, the `space of opens' is less so.
In the flawed argument for the existence of a space of opens from \autoref{sec:intro}, one assumes that given an open and a point, it is possible to semi-decide membership of the point in the open. But an open merely \emph{describes} the verifiable property under consideration --- it is a \emph{machine} that actually carries out this verification.

In the following subsection, we will show that the quotient map $\Sigma^{\Sigma^G} \twoheadrightarrow \O X$ has a section whenever the space $X$ is locally compact. From the above perspective, this fact perfectly explains why in this case a continuous evaluation map exists for $\O X$: given an open $u$ and a point $x$, one can apply the section to $u$ to acquire a machine $s(u)$ which can then verify membership of $x$ in $u$.

On the other hand, when $X$ is not locally compact the argument breaks down and we will see that there is no way to continuously map opens to machines in this case.

\begin{rem}\label{rem:nondet}
  Note that since we do not care which specific machine representing a given open is selected, we can relax the requirement of having an actual section to having a sufficiently well-behaved quotient. In particular, suplattice and preframe homomorphisms between frames are known to correspond to angelically and demonically nondeterministic maps respectively \cite{sondergaard1992non,winskel1983note,vickers2004universal} and it would seem that a nondeterministic choice of representative for each open would be sufficient to construct the evaluation map. Indeed, open and proper quotients, or more generally, triquotients, have such `non-deterministic sections' and these are strong enough conditions on the quotient map to entail that $X$ is exponentiable. We will discuss the role of nondeterminism in more detail in a later paper.
\end{rem}

\subsection{Relationship between machine space and the space of opens}\label{sec:machine_spaces_and_space_of_opens}

It is enlightening to consider how $\Sigma^{\Sigma^G}$ relates to $\Sigma^X$ in the case that the latter does exist. Intuitively, we expect there to be a quotient map $q\colon \Sigma^{\Sigma^G} \to \Sigma^X$ which sends a machine to its corresponding open. 

\begin{defi}\label{def:quotient}
 Let $X$ be a locally compact locale with $\O X = \langle G \mid R \rangle$. We may define the map $q$ of machine space via the universal property of the exponential $\Sigma^X$ applied to $\evweak$.
 
 \begin{center}
  \begin{tikzpicture}[node distance=3.5cm, auto]
    \node (A) {$\Sigma^{\Sigma^G} \times X$};
    \node (B) [below of=A] {$\Sigma^X \times X$};
    \node (C) [right =3.5cm of B] {$\Sigma$};
    \node (D) [left =1.7cm of B] {$\Sigma^X$};
    \node (E) [above of=D] {$\Sigma^{\Sigma^G}$};
    \node (F) [above of =C] {$\Sigma^{\Sigma^G} \times \Sigma^G$};
    \draw[->] (A) to node [swap] {$q \times X$} (B);
    \draw[->] (B) to node [swap] {$\ev$} (C);
    \draw[->] (A) to node {$\evweak$} (C);
    \draw[->] (F) to node {$\ev$} (C);
    \draw[->] (A) to node {$\Sigma^{\Sigma^G} \times i_X$} (F);
    \draw[->] (A) to node {$\evweak$} (C);
    \draw[dashed,->] (E) to node [swap] {$q$} (D);
  \end{tikzpicture}
\end{center}
Equivalently, $q$ is given by $\Sigma^{i_X}$ where $i_X$ is the inclusion of $X$ into $\Sigma^G$.
\end{defi}

This map behaves in accordance with our intuition. The diagram can be interpreted as meaning $\evweak(m,x) = \ev(q(m),x)$ which is to say that $m$ halts on $x$ if and only if $x$ lies in $q(m)$.
Moreover, it can be shown to be a quotient map. In fact, it necessarily has a section.
(Compare the $\Sigma$-split inclusions of Taylor \cite{taylor2002subspaces}.)

\begin{prop}\label{prop:section_from_exponential_to_machine_space}
Let $X$ be a locally compact locale with $\O X = \langle G \mid R \rangle$. The map $q \colon \Sigma^{\Sigma^G} \to \Sigma^X$ as defined in \autoref{def:quotient} has a section $s$ satisfying $\evweak \circ (s \times X) = \ev$.
\end{prop}

\begin{proof}
We define $s$ from the `weak' universal property of the weak exponential applied to $\ev$.

\begin{center}
  \begin{tikzpicture}[node distance=3.5cm, auto]
    \node (A) {$\Sigma^{\Sigma^G} \times X$};
    \node (B) [below of=A] {$\Sigma^X \times X$};
    \node (C) [right of=B] {$\Sigma$};
    \node (D) [left =1.7cm of B] {$\Sigma^X$};
    \node (E) [above of=D] {$\Sigma^{\Sigma^G}$};
    \draw[transform canvas={xshift=-0.75ex},->] (A) to node [swap,yshift=-7pt] {$q \times X$} (B);
    \draw[transform canvas={xshift=0.75ex},->] (B) to node [swap,yshift=-7pt] {$s \times X$} (A);
    \draw[->] (B) to node [swap] {$\ev$} (C);
    \draw[->] (A) to node {$\evweak$} (C);
    \draw[transform canvas={xshift=-0.75ex},dashed,->] (E) to node [swap] {$q$} (D);
    \draw[transform canvas={xshift=0.75ex},dashed,->] (D) to node [swap] {$s$} (E);
  \end{tikzpicture}
\end{center}

It is immediate that $\evweak \circ (s \times X) = \ev$ and so all that remains is to prove that $s$ is a section of $q$. 

By the diagram we have $\ev \circ (q \times X) \circ (s \times X) = \evweak \circ (s \times X) = \ev = \ev \circ (\id_{\Sigma^X} \times X)$. Hence by the uniqueness condition of the universal property of the `true' exponential $\Sigma^X$, we have $q \circ s = \id_{\Sigma^X}$, as required.
\end{proof}

The section $s$ continuously `picks out representatives of the equivalence classes defined by the quotient'. It can be viewed as sending an open $u$ to a machine which carries out the procedure of verification for $u$.

Conversely, if opens can be continuously assigned to machines, then we expect the exponential to exist. In fact, as we will see, it suffices for the map $q$ to be a `sufficiently good' quotient map.
(For example, we might only have a `non-deterministic section'.)

Recall that if a locale $X$ is given by $\O X = \langle G \mid R\rangle$ then the frame $\O X$ is a quotient of $\O(\Sigma^G)$ by the congruence generated from the relations $R$. Let us call this frame quotient map $\overline{q}$ for notational convenience, but as a function it is the same as $i_X^*$. Recall that $\O(\Sigma^G)$ equipped with the Scott topology is precisely the machine space $\Sigma^{\Sigma^G}$. Below we will write $\Sigma^{\Sigma^G}$ for $\O(\Sigma^G)$ equipped with the Scott topology even when viewing it as a topological space (as opposed to a locale). There is no danger in doing this since $\Sigma^{\Sigma^G}$ is spatial. We can now use $\overline{q}$ to induce a topology on $\O X$.
We will make use of the following result.
\begin{lem}
 Let $L$ and $M$ be posets admitting all joins and let $q\colon L \twoheadrightarrow M$ be a join-preserving surjection. Then if $L$ is equipped with the Scott-topology, the quotient topology induced on $M$ is also the Scott-topology.
\end{lem}
\begin{proof}
 We need to show that $U \subseteq M$ is Scott-open if and only if $q\inv(U)$ is Scott-open.
 If $U$ is Scott-open, then $q\inv(U)$ is Scott-open, because $q$ is Scott-continuous (as it preserves all joins).

 Now suppose $q\inv(U)$ is Scott-open. Consider $a \in U$ and suppose $a \le b$. Take $a' \in q\inv(\{a\})$ and $b' \in q\inv(\{b\})$. We may assume $b' \ge a'$ by joining $a'$ if necessary.
 Then $a' \in q\inv(U)$ and hence $b' \in q\inv(U)$, since $q\inv(U)$ is upward closed. Thus, $b = q(b') \in U$ and so $U$ is an upset.
 Next suppose $\dirsup D \in U$. Set $D' = q\inv({\downarrow} D)$. We note $D'$ is directed. Certainly $D'$ is nonempty. If $c,d \in D'$ then there is an $e \in D$ with $q(c), q(d) \le e$ by directedness of $D$. Thus, $q(c \vee d) = q(c) \vee q(d) \le e \in {\downarrow} D$ and so $c \vee d \in D'$.
 Furthermore, note that $q(D') = {\downarrow} D$ and $q(\dirsup D') = \dirsup q(D') = \dirsup {\downarrow} D = \dirsup D \in U$. Thus, $\dirsup D' \in q\inv(U)$ and so there is a $d' \in D'$ with $d' \in q\inv(U)$, since $q\inv(U)$ is Scott-open. Now $q(d') \in {\downarrow} D$ so that $q(d') \le d$ for some $d \in D$. Recalling that $U$ is an upset, we have $d \ge q(d') \in U$ and so $d$ is the desired $d \in U$. Thus, $U$ is Scott-open.
\end{proof}

\begin{cor}
 The quotient topology induced on $\O X$ by the function $\overline{q}\colon \Sigma^{\Sigma^G} \twoheadrightarrow \O X$ is the Scott topology on $\O X$.
\end{cor}

In the case that $X$ is locally compact, this map $\overline{q}$ agrees with the map $q$ we previously defined by the universal property.

\begin{lem}\label{lem:qbar_is_q_when_locally_compact}
The map $\overline{q}\colon \Sigma^{\Sigma^G} \to \O X$ is equal to $q\colon \Sigma^{\Sigma^G} \to \Sigma^X$ as defined in \autoref{def:quotient} whenever $X$ is locally compact.
\end{lem}

\begin{proof}
We first note that the codomains agree since $\O X = \Sigma^X$ for locally compact $X$.
Now recall that $q = \Sigma^{i_X}$. The natural transformation $(-)^{i_X}$ is the mate of $(-) \times i_X$ with respect to the adjunction $(-) \times X \dashv (-)^X$ and so we have the following commutative diagram.
\begin{center}
 \begin{tikzpicture}[node distance=3.5cm, auto]
  \node (A) {$\Hom(1, \Sigma^X)$};
  \node (B) [right of=A] {$\Hom(1 \times X, \Sigma)$};
  \node (C) [below=1.75cm of A] {$\Hom(1,\Sigma^{\Sigma^G})$};
  \node (D) [right of=C] {$\Hom(1 \times \Sigma^G, \Sigma)$};
  \node (E) [right of=B] {$\Hom(X, \Sigma)$};
  \node (F) [right of=D] {$\Hom(\Sigma^G, \Sigma)$};
  \draw[->] (A) to node {$\cong$} (B);
  \draw[->] (C) to node {$\Hom(1,\Sigma^{i_X})$} (A);
  \draw[->] (D) to node [swap] {$\Hom(1 \times i_X, \Sigma)$} (B);
  \draw[->] (C) to node {$\cong$} (D);
  \draw[->] (B) to node {$\cong$} (E);
  \draw[->] (D) to node {$\cong$} (F);
  \draw[->] (F) to node [swap] {$\Hom(i_X, \Sigma)$} (E);
 \end{tikzpicture}
\end{center}
The map $\Hom(i_X, \Sigma)$ corresponds to the action of the frame map $i_X^* = \overline{q}$, while $\Hom(1,\Sigma^{i_X})$ gives the action of $\Sigma^{i_X}$ on points.
Thus, we have that $\Sigma^{i_X}$ and $\overline{q}$ agree on points and hence they coincide.
\end{proof}

The failure of $\O X$ with the Scott topology to be an exponential $\Sigma^X$ when $X$ is not locally compact is due to this quotient map being badly behaved in general.
Indeed, if there were a section to $q$ as in  \autoref{prop:section_from_exponential_to_machine_space} we could show $X$ to be locally compact after all.

\begin{prop}\label{prop:section_to_machine_space_quotient_gives_exponential}
If $\overline{q}\colon \Sigma^{\Sigma^G} \to \O X$ has a continuous section $s \colon \O X \to \Sigma^{\Sigma^G}$, then $X$ is locally compact. Moreover, the evaluation map for $\O X$ is $\evweak \circ (s \times X)$.
\end{prop}

We thank the reviewer for suggesting the following proof, which is simpler than the one we originally had.
\begin{proof}
Since $\Sigma^G$ is locally compact, $\O(\Sigma^G)$ is a continuous frame and in particular a continuous dcpo. But continuous dcpos are known to be stable under Scott-continuous retracts and so $\O X$ is also a continuous dcpo (see \cite[Theorem 3.1.4.]{abramsky1994domain}).
Thus, $\O X$ is a continuous frame and $X$ is locally compact.

To see that $\evweak \circ (s \times X)$ is the evaluation map, first note that, by \autoref{lem:qbar_is_q_when_locally_compact}, $\overline{q} = q$ as defined in \autoref{def:quotient}. So by the commutative diagram in \autoref{def:quotient}, the evaluation map $\ev\colon \O X \times X \to \Sigma$ satisfies $\ev \circ (q \times X) = \evweak$.
Hence, $\evweak \circ (s \times X) = \ev \circ (q \times X) \circ (s \times X) = \ev$.
\end{proof}

In fact, we can say more (see \autoref{rem:nondet}). It is known that the lack of exponential objects is related to the fact that quotients are not stable under products (see \cite{day1970topological}). We have the following result. (This result requires us to assume $X$ is spatial and that we are in the category of topological spaces, though it is still useful for intuition. But see also \autoref{rem:exponential_section_presheaf}.)
\begin{prop}
 If $\overline{q} \times X\colon \Sigma^{\Sigma^G} \to \O X \times X$ is a quotient map, then $\O X$ (with the Scott topology) is the exponential $\Sigma^X$.
\end{prop}
\begin{proof}
 Consider the following diagram.
 \begin{center}
  \begin{tikzpicture}[node distance=4cm, auto]
    \node (D) {$\Sigma^{\Sigma^G} \times X$};
    \node (C) [left of=D] {$A \times X$};
    \node (E) [right of=D] {$\O X \times X$};
    \node (F) [below=2cm of D] {$\Sigma$};
    \draw[->] (E) to node {$\ev$} (F);
    \draw[->] (D) to node [swap] {$\evweak$} (F);
    \draw[->] (C) to node {$h' \times X$} (D);
    \draw[->] (D) to node {$\overline{q} \times X$} (E);
    \draw[->] (C) to node [swap] {$h$} (F);
    \draw[bend left,->] (C) to node {$h'' \times X$} (E);
  \end{tikzpicture}
 \end{center}
 Here $\ev$ is defined as the evaluation map from the exponential of the underlying sets. It is continuous since the composite $\evweak = \ev \circ (\overline{q} \times X)$ is continuous and $\overline{q} \times X$ is a quotient map.
 
 If $h \colon A \times X \to \Sigma$ then by the weak universal property of $\evweak$ we have a map $h'\colon A \to \Sigma^{\Sigma^G}$ making the diagram commute.
 Then $h'' = \overline{q}h'$ satisfies the condition needed for the weak universal property of $\ev\colon \O X \times X \to \Sigma$ to hold. Moreover, it is unique since we have uniqueness on the underlying sets by the universal property of the exponential in $\Set$. Thus, $\O X$ satisfies the universal property of the exponential $\Sigma^X$, as required.
\end{proof}

\begin{rem}\label{rem:exponential_section_presheaf}
  Another way to formulate these ideas is to embed $\Loc$ into the large presheaf category $\Set^{\Loc\op}$, where exponentials of locales exist. Then if the map $\Sigma^{\Sigma^G} \to \Sigma^X$ induced by $X \hookrightarrow \Sigma^G$ has a section, it can be shown that $X$ is locally compact. In fact, the analogue of a dcpo section (see \cite{vickers2004universal}) suffices. We omit the details.
\end{rem}

\section{Compactness and universal quantification}\label{sec:compactness}

We believe that machine space also has the potential to be useful in clarifying other aspects of topology. In this section we use it to study compactness. Escardó has explained compactness in terms of an algorithm for universal quantification in the setting of semantics of programming languages \cite{escardo2004synthetic,escardo2007infinite,escardo2008exhaustible}. Our approach yields a version of this algorithm that involves only topological/localic concepts and the core operations provided by the verifiability interpretation, without restricting us to a particular programming language.
Thus, we give a positive answer to the second question posed in the introduction: there is a uniform algorithm for universal quantification over any compact space.

\subsection{An algorithm for universal quantification}

Recall that a locally compact locale $X$ is compact precisely when $\{1\}$ is open in $\Sigma^X$. More generally, if $X$ is embedded into a locally compact space $Y$, then $X$ is compact if and only if there is an open of $\Sigma^Y$ consisting of the opens that cover $X$ (see \cite{hofmann1981local}).
In terms of machine space this says that if $\O X = \langle G \mid R \rangle$ then $X$ is compact if and only if there is an open in $\Sigma^{\Sigma^G}$ consisting of the machines which halt on all of $X$.
From the perspective of verifiability, this means we can semi-decide if a given machine always halts on $X$.
Indeed, we will provide an algorithm to carry out this very procedure.

It is important to understand exactly what it is such an algorithm needs to do.
Of course, to give such a procedure we need to know a precise description of the space $X$ by generators and relations. From this we can mathematically derive whether a particular formal combination of generators is equal to $1$.
On the other hand, recall that `the exact composition' of the machines in machine space is opaque to us. Indeed, the opens of machine space only allow us to test very particular properties of the machines. 
In the first case we are dealing with discrete `syntax', while in the second we are working with abstract machines which we do not have knowledge about a priori (see \autoref{rem:axiomatic_vs_verifiable}).

We are now ready to consider the algorithm for universal quantification over a compact locale $X$ with presentation $\langle G \mid R \rangle$. (To actually run such an algorithm we should restrict $G$ to be countable, but we can still imagine uncountable parallelism in theory.)

\begin{algorithm}
\caption{Semi-decision procedure for universal quantification over $X$}\label{alg:compactness}
\begin{algorithmic}
\State $\forall_X\colon \Sigma^{\Sigma^G} \to \Sigma$
\Function{$\forall_X$}{$m$}
    \For{each $S \in \powerfin(\powerfin(G))$}\Comment{performed in parallel}
        \If{$\bigvee_{F \in S}\bigwedge_{g \in F}g \sim 1$ with respect to $R$} \\ \vspace{-10pt}
            \For{each $F \in S$}
                \State $\mathtt{test}(m \in {\boxtimes} F)$\Comment{semi-decides if a branch of $m$ is contained in $F$}
            \EndFor
            \State $\mathtt{HALT}$
        \EndIf
    \EndFor
\EndFunction
\end{algorithmic}
\end{algorithm}

Here $\powerfin(T)$ denotes the set of finite subsets of $T$ and $\mathtt{test}(m \in {\boxtimes} F)$ halts precisely when $m \in {\boxtimes} F$. Recall that a machine $m = \bigvee_{i \in I} \bigwedge_{j\in J_i} g_j$ lies in ${\boxtimes} F$ if and only if there exists some $J_k$ such that $g_j \in F$ for each $j \in J_k$ --- that is, if every generator in some branch of $m$ lies in $F$.

Note that the algorithm itself involves testing whether $\bigvee_{F \in S}\bigwedge_{g \in F} g$ is a cover of $X$. This is on the `axiomatic' level of \autoref{rem:axiomatic_vs_verifiable} and how to test this is determined mathematically before running the algorithm. The role of the algorithm is to translate from the axiomatic level to the level of opaque machines.

\begin{thm}\label{thm:algorithm_correctness}
Let $X$ be a compact locale with $\O X = \langle G \mid R \rangle$.
Algorithm~\ref{alg:compactness} semi-decides if the given machine $m$ halts on all of $X$.
\end{thm}

\begin{proof}
Let $m = \bigvee_{i \in I} \bigwedge_{j\in J_i} g_j$. First suppose this covers $X$ so that $\bigvee_{i \in I} \bigwedge_{j\in J_i} g_j = 1$ in $X$. Since $X$ is compact there is a finite $I' \subseteq I$ such that $\bigvee_{i \in I'} \bigwedge_{j\in J_i} g_j = 1$. Now set $S = \{J_i \mid i\in I'\} \in \powerfin(\powerfin(G))$. This set $S$ will be considered in some parallel branch of the algorithm and we have $\bigvee_{F \in S}\bigwedge_{g \in F}g \sim 1$ in $R$ by construction.
Note that $\mathtt{test}(m \in {\boxtimes} J_i)$ halts for each $J_i \in S$, since $\bigwedge J_i \le m$ and hence $m \in {\boxtimes} J_i$. Therefore, this branch will reach $\mathtt{HALT}$ and so the entire computation halts.

Conversely, if the computation halts on $m$, then there is some $S$ which provides a finite refinement of the cover given by $m$ which covers $X$. Hence $\bigvee_{i \in I} \bigwedge_{j\in J_i} g_j$ covers $X$ and $m$ halts on all of $X$, as required.
\end{proof}

\begin{rem}
The usual way to link compactness to universal quantification is via the characterisation in terms of closed product projections (for example, see \cite{bauer2009dedekind,escardo2005notes,manuell2023pointfree}).
Thus, it is perhaps remarkable that this algorithm instead makes use of the standard open cover definition.
\end{rem}

It is worth discussing $\mathtt{test}(m \in {\boxtimes}F)$ in more detail.
Recall that in general an open merely suggests the existence of a procedure to semi-decide some membership and is not the procedure itself. Hence it is not a priori clear that there is some uniform way to compute $\mathtt{test}(m \in {\boxtimes}F)$. However, note that $\Sigma^{\Sigma^{G}}$ is locally compact and hence we expect there to be continuous assignment of opens to machines.

In Model~\ref{mod:bot} the procedure associated to ${\boxtimes} F$ is simply the act of observing the machine and comparing the generators it visits to the finite collection $F$.
It is clear that this can be done uniformly in $F$.
Alternatively, we might regard a machine as a program $m\colon \Sigma^G \to \Sigma$. To test whether a branch is contained in a subset $F$, we simply supply as input the generalised point $p_F\colon G \to \Sigma$ which halts precisely on the elements of $F$.

\begin{exa}[Cantor space]\label{ex:cantor}
Cantor space is the space $2^\N$ of infinite binary sequences. A presentation is given by $\langle z_n, u_n,\, n \in \N \mid z_n \wedge u_n = 0, z_n \vee u_n = 1 \rangle$. Intuitively, $z_n$ should be thought of as the open consisting of the sequences whose $n^\text{th}$ digit is $0$ and $u_n$ the open of sequences with $1$ in the $n^\text{th}$ position.
A machine will be given by a formal expression of the form $\bigvee_{i \in I} \left(\bigwedge_{j \in J_i} z_j \wedge \bigwedge_{k \in K_i} u_k\right)$.

Note that we can easily check if a \emph{finite} join $\bigvee_{i \in I} \left(\bigwedge_{j \in J_i} z_j \wedge \bigwedge_{k \in K_i} u_k\right)$ is a cover of $2^\N$ by distributing the joins over the meets and checking that each conjunct is a cover. A conjunct will be a cover if and only if it contains both $z_i$ and $u_i$ as disjuncts for some $i \in \N$.

The implementation of $\mathtt{test}(m \in {\boxtimes}F)$ depends on the precise model.
In any case, we arrive at an algorithm which achieves the same ends as that given by Escardó in \cite{escardo2004synthetic,escardo2008exhaustible}.
\end{exa}
We will explore how our approach relates to other representations of spaces and the link to Escardó's setting in \autoref{sec:dcpo}.

\begin{exa}[The closed interval]
 Consider the case of the closed real interval $X = [0,1]$. We can take a generating set consisting of certain intervals with rational end points, $G = \{[0,q) \mid q \in \Q \cap (0,1]\} \cup \{(q,1] \mid q \in \Q \cap [0,1)\}$. Then in the computational model, $[0,1]$ embeds as a subset of partial functions on $G$ and so a real $r$ corresponds to a function that halts on $[0,q)$ or $(q,1]$ if and only if $r \in [0,q)$ or $r \in (q,1]$ respectively. Then for a `machine' $m\colon \Sigma^G \to \Sigma$ the algorithm halts if $m$ halts on all of the (functions corresponding to) reals in $[0,1]$. It is elementary to compute when a formal join of rational intervals covers $[0,1]$.
 Finally, $\mathtt{test}(m \in {\boxtimes} F\})$ is given by simply calling $m(p_F)$ where $p_F \in \Sigma^G$ is given by $p_F(g) = \top \iff g \in F$.
\end{exa}

Another element to discuss further is the process of determining whether $\bigvee_{F \in S}\bigwedge_{g \in F}g \sim 1$ with respect to $R$.
We must do this by hand in order to construct the explicit algorithm for a given locale.
In concrete cases this is typically straightforward (as it is for the examples above), but ideally we might hope to have a general algorithm that does this for us. (In fact, it suffices for the algorithm to only semi-decide if the expression is a cover.) However, generating a congruence from given relations can involve a transfinite procedure and is hence, at the very least, a computationally nontrivial process.

Here it is natural to consider links to \emph{formal topology}, which is a predicative approach to topology that considers as its primary objects of study structures which correspond to presentations of locales. Of relevance here are the paper \cite{coquand2003inductively}, which discusses inductive generation of topologies, and \cite{vickers2006compactness}, which covers compactness and the relation to locale theory.
Theorem 15 of \cite{vickers2006compactness} gives conditions in terms of the generators for a set of formal finite joins of opens to be the set of all finite covers of a locale. Moreover, Proposition 11 of the same paper implies that the presentation of a compact locale may always be given in a way such that compactness is manifest and in which case it is easy to find the finite covers algorithmically.
In general, a proof of compactness is likely to lead to a description of the finite covers.

\begin{rem}\label{rem:overtness}
Constructively, there is a `dual' notion to compactness called \emph{overtness} (see \cite{taylor2011foundations,galoisTheoryGrothendieck})
which has the same relation to existential quantification as compactness has to universal quantification.
One can give an algorithm for existential quantification that is very similar to Algorithm~\ref{alg:compactness} except that instead of iterating over covers we iterate over `positive' (classically: non-zero) elements.
Classically, every locale is overt, though the resulting algorithm is still more subtle than naively searching through every real number in turn. Instead, it would use the fact that non-trivial open intervals with rational endpoints form a countable base for $\R$, which has the advantage of avoiding uncountable parallelism.

If the space $X$ is both compact and overt, we can combine the algorithms for universal and existential quantification to obtain an algorithm that determines whether a decidable predicate (i.e.\ clopen) holds everywhere on $X$ or not in finite time. This is analogous to the case that is of primary interest to Escardó.
\end{rem}

\subsection{Topological consequences}

It is possible to unwind Algorithm~\ref{alg:compactness} to give an open of $\Sigma^{\Sigma^G}$. Each expression $\mathtt{test}(m \in {\boxtimes}F)$ is replaced with the corresponding open ${\boxtimes}F$. The second \textbf{for} loop takes the meet of these opens, $\bigwedge_{F \in S} {\boxtimes}F$.
The first \textbf{for} (together with the \textbf{if}) corresponds to taking the join over all $S \in \powerfin(\powerfin(G))$ such that $\bigvee_{F \in S}\bigwedge_{g \in F}g$ covers $X$. Let us denote the set of such $S$ by $\mathcal{C}$. Thus, the resulting open is $\bigvee_{S \in \mathcal{C}}\bigwedge_{F\in S}{\boxtimes}F$. When $X$ is compact, this is the open of machines which cover $X$.
The following result is then immediate from \autoref{thm:algorithm_correctness}.

\begin{cor}\label{thm:open_of_covers}
Let $X$ be a compact locale with $\O X = \langle G \mid R \rangle$ and let $\mathcal{C}$ be the set $S \in \powerfin(\powerfin(G))$ satisfying $\bigvee_{F \in S}\bigwedge_{g \in F}g = 1$ in $X$. Then the open $\bigvee_{S \in \mathcal{C}} \bigwedge_{F\in S} {\boxtimes}F$ contains precisely the machines which cover $X$.
\end{cor}

Note that we can still define this open for non-compact $X$; however, in general it will only contain machines which correspond to covers which have a finite subcover.  

Also note \autoref{thm:open_of_covers} easily generalises to arbitrary compact sublocales $K$ of a general locale $X$ by simply viewing $K$ as direct sublocale of $\Sigma^G$.

\begin{rem}
 \autoref{thm:open_of_covers} also follows from the Hofmann--Mislove theorem \cite{hofmann1981local}.
 Moreover, the Hofmann--Mislove theorem also provides a converse giving that the set of all covering machines of $X$ is Scott-open if and only if $X$ is compact. (See also \cite[Lemma 7.4]{escardo2004synthetic}.)
 However, we believe restriction to machine space and the explicit algorithm provide an interesting new perspective on the situation.
\end{rem}

\begin{rem}
 Of course, in the overt case we can similarly extract an open of all `positive' machines for $X$ from the algorithm mentioned in \autoref{rem:overtness}.
\end{rem}

\subsection{Link to domain-theoretic approaches}\label{sec:dcpo}

Escardó approaches his compactness algorithm from the point of view of programming language semantics and domain theory. It is interesting to consider how our more topological approach compares to this one.
(Until now our pointfree results have been constructively valid. For simplicity, we will work classically in this section.)

Recall that data types in programming languages can be modelled as certain kinds of dcpos. The dcpos that arise in denotational semantics are usually \emph{continuous}, in which case the Scott topology is particularly well-behaved. We also note that the resulting topological spaces are locally compact.

Escardó considers the data type of infinite binary sequences. We imagine these are implemented as streams where the terms are computed in order one by one on demand.
Such data consists not only of infinite sequences, but also additional elements corresponding to sequences whose first $n$ terms are defined, but for which accessing later terms in the sequence causes the computation to diverge.

The elements of this data type form a dcpo ordered by definability so that, for instance, a sequence which hangs after producing $0,1,0$ is smaller than one which hangs after $0,1,0,0$. Equipping this with the Scott topology we obtain a space containing Cantor space as the subspace of completely defined sequences.

Our embedding of a locale $X$ into the space $\Sigma^G$ for a set of generators $G$ can be thought of as a reasonably canonical way to obtain additional `partially defined' elements in a similar way.
In fact, since the space $\Sigma^G$ is simply the power set $\powerset(G)$ (which is a continuous dcpo) equipped with the Scott topology, it can be understood as corresponding to the data type of partial functions from a discrete data type with $|G|$ elements to the unit type.

In this way our construction is an example of a \emph{domain embedding} where a topological space is represented as a subquotient of a dcpo.
In the domain theory literature, there is some desire for spaces to occur as the subspace of maximal elements of a dcpo, but this not the case for our embedding, and indeed, this is impossible if the space is not $T_1$. However, we can get closer to this ideal if we take the closure of $X$ as a subspace in $\powerset(G)$. A closed subset of a dcpo is a downward closed set which is closed under directed joins.
Moreover, a closed subspace of a continuous dcpo is a continuous dcpo and the Scott topology on the subset agrees with the subspace topology.
Thus, the closure of $X$ in $\powerset(G)$ is a dcpo which contains the points of $X$ and possibly some `more undefined' approximations to them and since this closure is the smallest Scott-closed set containing $X$, the points of $X$ lie as close to the top of the dcpo as possible in some sense.

\begin{exa}\label{ex:cantor_embedding}
 Let us consider the case of Cantor space from \autoref{ex:cantor}.
 Recall that $\O(2^\N)$ can be expressed as a quotient of $\O(\Sigma^G)$ where $G = \{z_n \mid n \in \N\} \sqcup \{u_n \mid n \in \N\}$ by the congruence $C = \langle z_n \wedge u_n = 0, z_n \vee u_n = 1,\, n \in \N \rangle$. The closure of such a sublocale corresponds to the congruence generated by the pairs $(0,c) \in C$. In this case we obtain the following presentation for the closure: $\langle z_n, u_n,\, n \in \N \mid z_n \wedge u_n = 0\rangle$.
 This describes a dcpo the elements of which can be thought of as \emph{partial} functions from $\N$ to $2 = \{0,1\}$. The points of Cantor space correspond to the total functions, which do appear as the maximal elements of this dcpo.
\end{exa}

In general, the compactness algorithm works equally well if we replace $\powerset(G)$ with such a continuous dcpo $D$, or indeed any locally compact sublocale of machine space which contains $X$. The algorithm takes an element $m$ of $\Sigma^D$ as input and simply applies Algorithm~\ref{alg:compactness} to $s(m)$ where $s$ is a section of the quotient $\Sigma^{\Sigma^G} \to \Sigma^D$ from \autoref{prop:section_from_exponential_to_machine_space}.

Note that the only place $s(m)$ appears in the algorithm is at the step where we run $\mathtt{test}(s(m) \in {\boxtimes}F)$. Thus, we may replace this step with any equivalent computation involving~$m$.

In this setting it is natural to consider an intuitive model of such a dcpo as a data type in an idealised programming language. If we view $\Sigma^{\Sigma^G}$ from this perspective recall that $\mathtt{test}(m \in {\boxtimes}F)$ is given by $m(p_F)$ where $p_F\colon G \to \Sigma$ halts precisely on the elements of $F$.
Other examples can be handled in a similar way.

\begin{continueexample}{ex:cantor_embedding}
 Let us return to the example of Cantor space. Using the dcpo $A$ of partial functions from $\N$ to $2$ instead of machine space we have a compactness algorithm where $\mathtt{test}(s(m) \in {\boxtimes}F)$ is implemented as follows:
 for each $E \in \powerfin(F)$ which does not contain both $u_n$ and $z_n$ for any $n$, we run $m(f_E)$ where
 \[
 f_E(n) = \begin{cases}
 0 & \text{if $z_n \in F$,} \\
 1 & \text{if $u_n \in F$,} \\
 \bot & \text{otherwise}
 \end{cases}
 \]
 and halt if any of these do.
 
 The inclusion $i\colon A \to \Sigma^G$ maps $f \in A$ to a function which halts on $z_n$ if $f(n) = 0$ and halts on $u_n$ if $f(n) = 1$. The quotient map $q\colon \Sigma^{\Sigma^G} \to \Sigma^A$ sends $m$ to a map $x \mapsto m(i(x))$. 
 Note that for a given $m \in \Sigma^A$, the above procedure $\mathtt{test}(s(m) \in {\boxtimes}F)$ for $F \in \powerfin(G)$ can be easily extended to infinite $F \in \Sigma^G$. This gives the associated section $s\colon \Sigma^A \to \Sigma^{\Sigma^G}$ of~$q$.
\end{continueexample}

The above representation of Cantor space is not quite the one used by Escardó, since he views Cantor space not as a function, but as a \emph{stream} --- that is, there is an ordering on indices so that if the value at $i$ is defined, so are all previous values.
This is related to an alternative presentation of Cantor space in terms of prefixes.
\begin{exa}
 A presentation for $2^\N$ can be given with generators $\ell_p$ corresponding to sequences with finite non-empty prefix $p \in 2^+ \coloneqq 2^* \setminus \{\epsilon\}$. We have
\begin{align*}
 \O 2^\N = \langle \ell_p,\, p \in 2^+ \mid {} &\ell_p \wedge \ell_q = \ell_q \text{ for $p \prec q$}, \\
 & \ell_p \wedge \ell_q = 0 \text{ if $q \nprec p$ and $p \nprec q$}, \\
 & \ell_{p \# 0} \vee \ell_{p \# 1} = \ell_p \text{ for all $p$}, \\
 & \ell_0 \vee \ell_1 = 1 \rangle
\end{align*}
where $p \prec q$ means $p$ is a prefix of $q$ and $\#$ denotes concatenation.

However, taking the closure of this sublocale of $\Sigma^{2^+}$ does not yield the correct dcpo. This is because it ignores the natural order structure on the prefixes.
The solution is to modify the parent space $\Sigma^G$, so that instead of taking the free frame on a \emph{set} $G$, we take the free frame on a poset. In this case, we equip $2^+$ with the reverse of the prefix order.
Then $\O 2^\N$ is the quotient of this frame by the congruence generated by the second, third and fourth relations above (since the first relation is already handled by the order structure on $G$). The closure of this sublocale has the presentation
$\langle \{\ell_p : p \in (2^+)\op\} \in \Pos \mid \ell_p \wedge \ell_q = 0 \text{ if $q \nprec p$ and $p \nprec q$} \rangle$. It is not hard to see that this is the locale of finite and infinite binary sequences. (Indeed, it is the final two relations in the presentation for $2^\N$ above that forces all sequences to be infinite.)
This is precisely the dcpo considered by Escardó. (Finite sequences are interpreted as sequences which are undefined after some point.)

We can now see how our algorithm reduces to a (non-optimised) version of Escardó's algorithm.

Let us first consider how to check if a finite join $\bigvee_{i \in I} \left(\bigwedge_{p \in J_i} \ell_p \right)$ is a cover of $2^\N$.
We first compute the finite meets using the first two relations to obtain a join of basic generators (omitting zeros from the join).
Now let $N$ be the length of the longest prefix.
We pad each $p$ in the join with $0$s and $1$s to obtain all possible strings of length $N$ with prefix $p$. The join is a cover if and only if the resulting set contains every string of length $N$.

Finally, $\mathtt{test}(s(m) \in {\boxtimes} F)$ is computed as follows. If $\bigwedge F = 0$ there is nothing to do. Otherwise, $\bigwedge F = \ell_p$ and we run $m$ on the finite sequence $p$.
\end{exa}

\begin{rem}
 The locale corresponding to the free frame on a poset $P$ is $\Sigma^{\mathscr{I}P}$ where $\mathscr{I}P$ is the \emph{ideal completion} of $P$. In the case above, the ideal completion of $(2^*)\op$ is actually already the desired dcpo. Nonetheless, it is still interesting to see that this fits into our general approach.
\end{rem}

\section*{Acknowledgements}

We thank the referees for their comments, which improved the presentation of this paper.

\bibliographystyle{alphaurl}
\bibliography{bibliography}
\end{document}